\documentclass[11pt]{article}
\usepackage{amssymb}

\newtheorem{thm}     {Theorem}[section]
\newtheorem{prop}    [thm]{Proposition}

\newtheorem{lemma}   [thm]{Lemma}

\newcommand{\proof} {\noindent{\bf Proof. }}

\newcommand{\B}{\mathbb B}
\newcommand{\C}{\mathbb C}
\newcommand{\D}{\mathbb D}

\newcommand{\R}{\mathbb R}

\begin{document}

\title{ On holomorphic mappings  of  strictly pseudoconvex domains}
\author{Alexandre Sukhov }
\date{}
\maketitle

\bigskip
 
{\small Abstract.  We study the boundary regularity of proper holomorphic mappings between strictly pseudoconvex domains with boundaries of class $C^2$.}

MSC: 32H02.

Key words: strictly pseudoconvex domain, totally real manifold, wedge,  proper holomorphic mapping, boundary regularity, complex disc.

\bigskip

{\small

 Universit\'e  de Lille, Laboratoire
Paul Painlev\'e,
U.F.R. de
Math\'ematiques, 59655 Villeneuve d'Ascq, Cedex, France, e-mail: sukhov@math.univ-lille1.fr,
 and  Institut of Mathematics with Computing Centre , Ufa Federal Research Centre of Russian
Academy of Sciences, 450077, Chernyshevsky Str. 112, Ufa, Russia.

The author is partially suported by Labex CEMPI.
}

\section{Introduction}

At present the old problem of boundary regularity of proper or biholomorphic mappings between strictly pseudoconvex domains  in $\C^n$ is almost completely clarified due to the  contribution of several authors. Theorem of Ch. Fefferman \cite{F}  asserts that a biholomorphic mapping betweeen strictly pseudoconvex domains with boundaries of class $C^\infty$ extends as a $C^\infty$ diffeomorphism between their closures. The original  proof is based on the study of asymptotic behavior of the Bergman kernel near the boundary and is difficult to adapt to the case of finite smoothness. Later several authors developed different approaches.  One of them is due to L. NIrenberg, S.Webster, P. Yang \cite{NWY} and uses a smooth version of the reflection principle. Pushing further their approach,  S.Pinchuk and S.Khasanov \cite{PiKh} proved that a proper holomorphic mapping between strictly pseudoconvex domains with $C^s$ boundaries, with $s > 2$  real, extends to the boundary as a mapping of class $C^{s-1}$ if $s$ is not an integer, and of class $C^{s-1-\varepsilon}$, with any $\varepsilon > 0$, when $s$ is an integer. Similar results are obtained by L. Lempert \cite{Le1,Le2} by quite different techniques of extremal discs for the Kobayashi metric. Y. Khurumov \cite{Kh} proved that an even a better  result still remains true   with the loss of regularity on $1/2$. 
However, a natural question on the precise regularity in the case where  the boundaries of class exactly $C^2$ (i.e. the minimal possible regularity of boundaries), remains open. Y.Khurumov announced without further details that his result remains true also in this case, but, to the best of my knowledge, a detailed  proof is not available. The well-known result  in this case was obtained independently by  G.Henkin \cite{Khen} and  S.Pinchuk \cite{Pi1}) and states that a mapping extends to the boundary as a H\"older $1/2$-continuous mapping. 




Our  main result is the following

\begin{thm}
\label{MainTheo1}
 Let $f:\Omega_1 \to \Omega_2$ be a proper holomorphic mapping between  bounded strictly pseudoconvex domains  $\Omega_j \subset \C^n$, $j = 1,2$,  with   boundaries  $b\Omega_j$ of class $C^2$. Then $f$ extends to a mapping of the H\"older class $C^{\alpha}(\overline{\Omega}_1)$
for each $\alpha \in [0,1[$.
\end{thm}
In particular, this result means that the above mentioned theorem of Lempert-Pinchuk-Khasanov still remains true for boundaries of class $C^2$. This result is obtained in \cite{Su1} under the additional assumption that the boundary $b\Omega_1$ is of class $C^{2 + \varepsilon}$, $\varepsilon > 0$. From this point of view the present paper is a continuation and the complement of \cite{Su1}. The proof consists of two main parts. The first one is based on estimates of the Kobayashi metric in a tube 
neighborhood of a totally real manifold (see \cite{ChSu}). These estimates allow to obtain uniform  H\"older estimates on analytic discs glued to a prescribed totally real manifold which in turn allows to obtain the H\"older regularity of holomorphic mappings between wedge-type domains with totally real edges.
This argument is presented in details in \cite{Su1}. The present paper is devoted to the second main step of the proof. We study geometric properties of analytic discs attached  along an arc to a totally real manifold of class $C^1$. This  {\it gluing disc argument}, which  is often quite useful for the study of totally real submanifolds,  was introduced in \cite{Pi2}.  In the case where the regularity of a totally real manifold is higher than $C^1$, this construction was developped by several authors. The $C^1$ case was  considered by E.Chirka \cite{KhenCh} and Y.Khurumov \cite{Kh2}. However, for our goals we need some additional properties of such discs which are not explicitely stated there. This is the reason why  I present some details for completeness of exposition. I stress that I do not claim any originality here. Only applications of this construction are new.

\section{Terminology and notations}

We briefly recall some well-known definitions and basic notations.

 Let $\Omega$ be a domain in $\C^n$. For a positive integer $k$, denote by  $C^k(\Omega)$  the space of $C^k$-smooth complex-valued functions in~$\Omega$. Also $C^k(\overline\Omega)$ denotes the class of functions whose partial derivatives up to order $k$ extend as continuous functions on 
$\overline\Omega$. Let  $s > 0$ be a real noninteger and let $k$ be its integer part. Then  $C^s(\Omega)$ denotes the space of functions of class $C^k(\overline\Omega)$ such that their partial derivatives of order $k$ are (globally)  $(s-k)$-H\"older continuous in $\Omega$; these derivatives  automatically satisfy the  $(s-k)$-H\"older condition on 
$\overline\Omega$ so the notation $C^s(\overline\Omega)$ for the same space of functions is also appropriate. We also use the space $L^p(\Omega)$ of Lebesgue $p$-integrable functions with the usual norm 
$\parallel f \parallel_{L^p(\Omega)}$; if $f = (f_1,...,f_m)$ is a vector valued function, we set $\parallel f \parallel_{L^p(\Omega)} = \sum_{j=1}^m \parallel f_j \parallel_{L^p(\Omega)}$.
We also denote by $W^{k,p}(\Omega)$  (resp. $W^{k,p}_{loc}(\Omega)$) the Sobolev spaces of (vector) functions with generalized derivatives of order up to $k$ which are Lebesgue $p$-integrable (resp. locally) in $\Omega$. 

Denote by $\D = \{ \zeta \in \C: \vert \zeta \vert < 1 \}$ the unit disc of $\C$. Recall that in this special case the Sobolev embedding theorem 
asserts that for $\alpha = 1 - 2/p$ and $p > 2$ the natural inclusion  $W^{1,p}(\D) \to C^{\alpha}(\D)$  is a linear bounded compact operator.


A (closed) real submanifold $E$ of a domain $\Omega \subset \C^n$ is of class $C^s$ (with real $s  \ge 1$) if for every point $p \in E$ there exists an open neighbourhood $U$ of $p$ and a map $\rho: U \longrightarrow \R^d$ of the maximal rank $d<2n$ and of class $C^s$  such that $E \cap U = \rho^{-1}(0)$; then $\rho$ is called a local defining (vector-valued ) function of $E$. The positive integer $d$ is the real codimension of $E$. In the most important special case $d=1$ we obtain the class of real hypersurfaces.

Let $J$ denotes the standard complex structure of $\C^n$. In other words, $J$ acts on a vector $v$ by multiplication by $i$ that is $J v = i v$. For every $p \in E$ the {\it holomorphic tangent space}  $H_pE:= T_pE \cap J(T_pE)$ is the maximal complex subspace of the tangent space $T_pE$ of $E$ at $p$. Clearly $H_pE =   \{ v \in \C^n: \partial \rho(p) v = 0 \}$. The complex dimension of $H_pE$ is called the CR dimension of $E$ at $p$; a manifold $E$ is called a {\it CR  (Cauchy-Riemann) manifold} if its CR dimension is independent of $p \in E$.

A real submanifold $E \subset \Omega$ is called {\it generic} (or generating) if the complex span of $T_pE$ coincides with $\C^n$ for all $p \in E$. Note that every generic manifold of real codimension $d$ is a CR manifold of CR dimension $n-d$. A function $\rho = (\rho_1,...,\rho_d)$ defines a generic manifold if $\partial\rho_1 \wedge ...\wedge \rho_d \neq 0$.
Of special importance are the so-called {\it totally real manifolds}, i.e., submanifolds $E$ for which $H_pE = \{ 0 \}$ at every $p \in E$. A totally real manifold in $\C^n$  is generic if and only if its real dimension is equal to $n$; this is the maximal possible value for the dimension of a totally real manifold. 

 Let $\Omega$ be a bounded domain in~$\C^n$. Suppose that its boundary $b\Omega$ is a (compact) real hypersurface of class $C^s$ in $\C^n$. Then there exists a $C^s$-smooth real 
function $\rho$ in a neighbourhood $U$ of the closure $\overline\Omega$ such that $\Omega = \{ \rho < 0 \}$ and 
$d\rho|_{b\Omega} \ne 0$. We call such a function $\rho$ a global defining function. If $s \ge 2$ one may consider 
{\it the Levi form} of $\rho$:
\begin{eqnarray}
\label{Leviform}L(\rho,p,v) = \sum_{j,k = 1}^n \frac{\partial^2\rho}{\partial z_j\partial\overline{z}_k}(p)v_j \overline v_k.
\end{eqnarray}
A bounded domain $\Omega$ with $C^2$ boundary is called  {\it strictly pseudoconvex}
if $L(\rho,p,v) > 0$  for every nonzero  vector $v \in H_p(b\Omega)$.

Consider a wedge-type domain 

\begin{eqnarray}
\label{wedge}
W = \{ z \in \C^n: \phi_j(z) < 0, j = 1,...,n \}
\end{eqnarray}
with the edge 

\begin{eqnarray}
\label{edge}
E = \{ z \in \C^n: \phi_j(z) = 0, j = 1,...,n \}
\end{eqnarray}
We assume that the defining functions $\phi_j$ are of class $C^{1}$. Furthermore, as usual we suppose that $E$ is a generic manifold that is 
$\partial \phi_1 \wedge ... \wedge \partial \phi_n \neq 0$ in a neighborhood of $E$. 

Given $\delta > 0$ (which is supposed to be small enough) we also define a shrinked wedge 
\begin{eqnarray}
\label{shrinked}
W_{\delta} = \{ z \in \C^n: \phi_j - \delta \sum_{l \neq j} \phi_l < 0, j= 1,...,n \} \subset W
\end{eqnarray}
It has the same edge $E$. Note that there exists a constant $C > 0$ such that for every point $z \in W_\delta$ one has 
\begin{eqnarray}
\label{dist}
C^{-1} dist(z, b W) \le dist(z, E) \le C dist (z,b W)
\end{eqnarray}

\section{Gluing complex discs to $C^1$ totally real manifolds} 
In this section we present the main technical tool for the proof of Theorem \ref{MainTheo1}. 
Consider a wedge-type domain (\ref{wedge})  with the edge (\ref{edge}).

We need the well-known construction of filling a wedge $W$ (or, more precisely, the wedge $W_\delta$) by complex discs gluing to $E$ along an open arc.  

A complex (or analytic, or holomorphic) disc is a holomorphic map $h: \D \to \C^n$ which is at least continous on the closed disc $\overline\D$. 
We say that such a disc is glued (or attached) to a subset $K$ of $\C^n$ along an (open, nonempty) arc $\gamma \subset b\D)$, if $f(\gamma) \subset K$.

Our presentation consists of   several steps.

\subsection{The generalized Bishop equation, the existence  and the regularity of discs}



Let $E$ be an $n$-dimensional totally real manifold of class $C^1$ in a neighborhood of $0$ in $\C^n$; we assume $0 \in E$. After a linear change of coordinates, using the implicit function theorem  we also may assume that in a neighbourhood $\Omega$ of the origin the manifold $E$ is defined by the (vector) equation 

\begin{eqnarray}
\label{edge2}
y = h(x)
\end{eqnarray}
 where a  vector function $h = (h_1,...,h_n)$ of class $C^1$ in a neighborhood of $0$ in $\R^n$ and satisfies the conditions 
 
 \begin{eqnarray}
 \label{edge3}
 h_j(0) = 0, \,\,\,\nabla h_j(0) = 0 , \, \, j = 1,...,n. 
\end{eqnarray}
Here and below $\nabla$ denotes the gradient.  

\bigskip


Fix a positive noninteger $s$. Consider  the Hilbert transform $T: u \to Tu$, associating to  a real function $u \in C^s(b\D)$  its harmonic conjugate function vanishing at the origin. In orther words, $u + iTu$ is a trace on $b\D$ of a function, holomorphic on $\D$ and of class $C^s(\D)$, and $Tu(0) = 0$. 

Recall that explicitely the Hilbert transform  is given  by 

\begin{eqnarray*}
Tu (e^{i\theta}) = \frac{1}{2\pi} v. p. \int_{-\pi}^{\pi} u(e^{it}) \cot \left ( \frac{\theta - t}{2} \right ) dt
\end{eqnarray*}

 This is a classical linear singular integral operator; it is  bounded on the space $C^s(b\D)$ for any non-integer $s > 0$. Furthermore, for $p > 1$ the operator $T:L^p(b\D) \to \L^p(b\D)$ is a bounded linear operator as well; we denote by $\parallel T \parallel_p$ its norm.

 Let $b\D^+ = \{ e^{i\theta}: \theta \in [0,\pi] \}$ and $b\D^- = \{ e^{i\theta}: \theta \in ]\pi, 2 \pi[ \}$ denote the upper and the lower semicircles respectively. Fix a $C^\infty$-smooth real functions $\psi_j$ on $b\D$ such that $\psi_j\vert b\D^+ = 0$ and  $\psi_j\vert b\D^- < 0$, $j=1,...,n$ (one may take the same function independent of $j$). Set $\psi = (\psi_1,...,\psi_n)$. Consider the {\it generalized Bishop equation} 
\begin{eqnarray}
\label{BishopEq}
u(\zeta) =  -Th(u(\zeta)) - tT\psi(\zeta) + c, \,\, \zeta \in b\D ,
\end{eqnarray}
where $c \in \R^n$ and $t = (t_1,...,t_n) \in \R^n$, $t_j \ge 0$, are real parameters. We will prove that for any $p > 2$, and for any $c$, $t$ close enough to the origin, this singular integral equation admits a unique solution $u(c,t)(\zeta)$ in the  Sobolev class $W^{1,p}(b\D)$. Such a solution is  of class $C^{\alpha}(b\D)$, $\alpha = 1 - 2/p$,  by the Sobolev embedding. 

\bigskip

Before proceed the proof, we explain how such a solution is related to complex discs glued to $E$ along $b\D^+$.
Indeed, consider the function $$U(c,t)\zeta) = u(c,t)(\zeta) + i h(u(c,t)(\zeta)) + i t\psi(\zeta).$$ 
Since $T^2 = -Id$ and $u$ is a solution of (\ref{BishopEq}), the function $U$ extends holomorphically on $\D$ as a function 
\begin{eqnarray}
\label{Discs}
H(c,t)(\zeta) = P U(c,t), \,\,\zeta \in \D
\end{eqnarray}
 of class $C^\alpha(\D)$.  Here $P$ denotes the Poisson operator of harmonic extension to $\D$:

\begin{eqnarray}
\label{Poisson}
PU(c,t)(\zeta) =  \frac{1}{2\pi} \int_{-\pi}^{\pi} \frac{1- \vert \zeta \vert^2}{\vert e^{it} - \zeta \vert^2} U(c,t)(e^{it})dt
\end{eqnarray}

The function  $\psi$ vanishes on $b\D^+$, so by (\ref{edge2}) we have $H(c,t)(b\D^+) \subset E$ for all $(c,t,d)$. 

\bigskip 

It is convenient to  extend the equation (\ref{BishopEq}) on the whole space. Fix a $C^\infty$ smooth 
function $\lambda: \R^{n} \to \R^+$ equal to $1$ on the unit ball $\B^{n}$ and vanishing on 
$\R^{n} \setminus 2\B^{n}$. For $\delta > 0$ small enough  the function 
$h_\delta(x) = \lambda(x/\delta)h(x)$  naturally extends (by $0$) on the whole $\R^{n}$.
Fix $\tau > 0$ small enough (it will be choosen later). Then in view of (\ref{edge3}) we can choose $\delta = \delta(\tau) > 0$ 
such that the gradient $\nabla h_\delta(x)$ is small on the whole $\R^{n}$:

\begin{eqnarray}
\label{edge4}
\parallel \nabla h_\delta \parallel_{L^\infty(\R^{n})} \le \tau
\end{eqnarray}

First we study the global equation 

\begin{eqnarray}
\label{BishopEq2}
u(\zeta) =    -Th_\delta(u(\zeta)) - tT\psi(\zeta) + c, \,\, \zeta \in b\D ,
\end{eqnarray}
We prove that its solutions depend continuously on parametrs $(c,t)$; this allows to localize  the  solutions and  to conclude with the 
initial equation (\ref{BishopEq}).

\bigskip

Let $V$ be a domain in $\R^m$ and $f \in L^p(V \times b\D)$. Then by the Fubini theorem $Tf \in 
L^p(V \times b\D)$ (the variables in $V$ are viewed as parameters when the operator $T$ acts). Hence, keeping the same notation, we obtain
a bounded linear operator $T: L^p(V \times b\D) \to L^p(V \times b\D)$  with the same norm as in $L^p(b\D)$.
We again denote its norm by $\parallel T \parallel_p$. 

Fix a domain   $V\subset \R^{2n}$  of the parameters $(c,t)$.

\begin{lemma}
\label{DiscExistence}
Under the above assumptions, for any $p > 1$,  one can choose $\tau > 0$ in (\ref{edge4}), and $\delta = \delta(\tau) > 0$,  such that   the equation (\ref{BishopEq2}) admits a unique solution $u(c,t)(\zeta) \in L^p(V  \times b\D)$.
\end{lemma}
\proof  Consider the operator $\Phi: L^p(V \times  b\D) \to L^p(V \times  b\D)$ defined by 
$$\Phi: u \mapsto  -Th_\delta(u(\zeta)) - tT\psi(\zeta) + c.$$ In view of (\ref{edge4}) it follows by the finite differences 
theorem that for all $u^1$, $u^2$ from $L^p(V \times  b\D)$ one has 
\begin{eqnarray*}
\parallel \Phi(u^1) - \Phi^2(u^2) \parallel_{L^p(V  \times b\D)} \le \parallel T \parallel_p \parallel h_\delta(u^1) - h_\delta(u^2) \parallel_{L^p(V  \times b\D)} \le (1/2) \parallel u^1 - u^2 \parallel_{L^p(V \times  b\D)}
\end{eqnarray*}
when $\tau$ in (\ref{edge4}) is fixed small enough.  Hence  $\Phi$ is a contracting map and lemma is proved.
\bigskip

Since $V$ is arbitrary we conclude that the equation (\ref{BishopEq2}) admits a unique solution $u \in L^p_{loc}(\R^{2n}  \times b\D)$. By this 
space we mean the space of $L^p$ functions on $K \times b\D$ for each (Lebesgue) measurable compact subset $K \subset \R^{2n}$.

\bigskip

Next we study the regularity of solutions of (\ref{BishopEq2}) in the Sobolev scale.  Let $\Omega$ be a domain in $\R^k$ and let  $f $ be a function in $\Omega$. 
Denote by $e_j$, $j=1,...,k$ the canonical basis of $\R^k$.
Given $j = 1,...,k$  and $\Delta x_j \in \R^*$ consider the finte differences 
$$\Delta f / \Delta x_j = (f(x + e_j \Delta x_j) - f(x))/\Delta x_j.$$
Recall two well-known properties of the Sobolev spaces (see, for example, \cite{Ni}).

\begin{itemize}
\item[(1)] If a function $f$ is in  $W^{1,p}_{loc}(\Omega)$, then for every subdomain $\tilde\Omega \subset \Omega$, the finite differences $\Delta f/\Delta x_j$ converge 
in $L^p(\tilde \Omega)$ as $\Delta x_j \to 0$ to the generalized derivative $\partial f/\partial x_j$.
\item[(2)] $f \in W^{1,p}_{loc}(\Omega)$ if and only if the finite differences $\Delta f/\Delta x_j$ are bounded 
in $L^p(\tilde \Omega)$ uniformly in $\Delta x_j$.
\end{itemize}

\begin{lemma}
\label{DiscRegularity}
Every solution $u$ of (\ref{BishopEq2}) is of class $W^{1,p}_{loc}(\R^{2n} \times b\D)$. 
\end{lemma}
\proof  Let $x_j$ denote one of the variables $c_j$, $t_j$  or $\zeta \in b\D$. We estimate the finite difference $\Delta u / \Delta x_j$. It follows from 
(\ref{edge4}) that $h_\delta$ satisfies the Lispchitz condition  with the Lipschitz constant $\tau$. Hence from (\ref{BishopEq2}) we have:
$$\parallel\Delta u/ \Delta x_j \parallel_{L^p(V \times  b\D)}\le C_1 \tau \parallel \Delta u/\Delta/x_j\parallel_{L^p(V \times b\D)} + C_2$$
where $C_j > 0$ are constants. When $\tau > 0$ (and so $\delta > 0$) are small enough, we obtain
$$\parallel\Delta u/\Delta x_j\parallel_{L^p(V \times  b\D)} \le C_3$$
for some constant $C_3 > 0$, and therefore $u \in W^{1,p}_{loc}(\R^{2n} \times b\D)$. Lemma is proved.
\bigskip

It follows by the Sobolev embedding that a solution $u$ belongs to $C^{1-(2n+1)/p}(V\times b\D)$, where $V$ is an open subset in $\R^{2n}$.
In particular, the constructed family of discs 
is continuous in all variables for $p$ big enough. Now we note that for $t = 0$ the equation (\ref{BishopEq2}) admits a constant solution $u(c,0)(\zeta) = c$. When $c$ is close enough to 
the origin in $\R^n$, this solution gives a point $c + ih(c) \in E$. By continuity and uniqueness of solutions, there exists a neighborhood $V$ of the origin in $\R^{2n}$, such that for  $(c,t) \in V$ any solution of (\ref{BishopEq2}) is a solution of (\ref{BishopEq}). We have proved the following

\begin{lemma}
\label{BishopLemma1}
Given $p > 2$ the exists a neighborhood $V$ of the origin in $\R^{2n} $ such that the Bishop equation (\ref{BishopEq}) admits a unique solution 
$u(c,t)(\zeta) \in W^{1,p}(V \times  b\D)$.
\end{lemma}
Since $p$ is arbitrary, we obtain that our equation admits solutions in the Holder class $C^{\alpha}(V \times  b\D)$ with $\alpha = 1-(2n+1)/p$. Note that here $V$ depends on $p$ (and hence, on $\alpha$). Nevertheless, it follows from \cite{ChSu} that for each  $(c,t)$ fixed, the map $\zeta \mapsto u(c,t)(\zeta)$ is of class $C^\alpha(b\D)$ for every $\alpha < 1$.

\subsection{Stability of discs} Until now we did not study any geometric properties of the family (\ref{Discs}). Here we consider some of them which will be useful for our applications.







We represent the family (\ref{Discs}) as a small perturbation in the $W^{1,p}$ norm of some model family.
The model case arises when $E = \R^n$ that is $h = 0$ in (\ref{edge}). Then the equation (\ref{BishopEq}) takes the from 

\begin{eqnarray}
\label{BishopEq7}
u(\zeta) =   - tT\psi(\zeta) + c, \,\, \zeta \in b\D ,
\end{eqnarray}
where as usual $c \in \R^n$ and $t = (t_1,...,t_n)$, $t_j \ge 0$, are real parameters. This equation is already solved with respect to $u$; one can view (\ref{BishopEq7}) as its general solution. In this case the family (\ref{Discs}) becomes 

\begin{eqnarray}
\label{Discs2}
H(c,t)(\zeta) = P U(c,t), \,\,\zeta \in \D
\end{eqnarray}
where

\begin{eqnarray}
\label{FlatDiscs}
U(c,t)\zeta) = - tT\psi(\zeta) + c +  i t\psi(\zeta).
\end{eqnarray}
This family arises  from the family of complex lines intersecting $\R^n$ along 
real lines; it is simply obtained by a biholomorphic reparametrization of arising half-lines by the unit disc. These lines are given by $l(c,t): \zeta \mapsto t\zeta + c$, $\zeta \in \C$. The conformal map $-T \psi + i\psi$ takes the unit disc into a smoothly bounded 
domain in the lower half-plane, gluing $b\D^+$ to the real axes. In what follows  we refer this case as the flat case.

When $h$ is of class $C^s$ with $s > 1$,  it follows easily by the implicit function theorem that the general family (\ref{Discs}) is obtained from the flat family (\ref{FlatDiscs}) by a small $C^s$-perturbation and hence has the similar geometric properties. The $C^1$ case requires some additional technical analysis.

 Let  $E$ be a totally real manifold given by (\ref{edge2}), (\ref{edge3}).  Given $d \in I \setminus \{ 0 \}$, where $I \ni 0$ is an open interval  small  enough, consider the manifolds $E_d$ given by
 
\begin{eqnarray}
\label{edge.d}
y = d^{-1}h(d x)
\end{eqnarray}
Note that for every $d  \neq 0$ the manifold $E_d$ is biholomorphic to $E$ via the isotropic dilation $z \mapsto  d^{-1}z$. Set $h(x,d) = d^{-1}h(dx)$. 
Thus, we consider the $1$-parameter family $E_d$ of totally real manifolds defined by the equation 
\begin{eqnarray}
\label{edge2.d}
y = h(x,d)
\end{eqnarray}
 

 \begin{eqnarray}
 \label{edge3.d}
 h_j(0,d) = 0, \,\,\,  \nabla_x h_j(0,d) = 0,  \, \, \ d \in I,  \, \, j = 1,...,n. 
 \end{eqnarray}
(we consider the gradient $\nabla_x$ with respect to $x$). 
 Hence for each $(c,t,d)$ we have the discs
$H(c,t,d)$  defined by (\ref{Discs}). If $d = 0$, we have $E_0 = \{ y = 0 \} = \R^n = T_0(E)$ that is, the flat case. By the uniqueness of the solution of the Bishop equation, 
the family $H(c,t,0)( \zeta)$ coincides with the family ( \ref{FlatDiscs}).

 \begin{lemma}
  \label{LemDeformation}
 For any $p > 1$  one has 
  $$ \parallel H(c,t,d)( \zeta) - H(c,t,0)( \zeta) \parallel_{W^{1,p}(V  \times  \D)}  \to 0$$
  as $d  \to 0$.
   \end{lemma}
 \proof Let $u^0$ be a solution of the flat Bishop equation ( \ref{BishopEq7}). For any $d$ let $u(c,t,d)$ be a solution of the Bishop equation 
  \begin{eqnarray}
  \label{BishopEq8}
 u = -T h(u,d) - tT \psi + c
  \end{eqnarray}
 It follows from (\ref{edge3.d}) that previous estimates on the norm of $u$ are uniform in $d$ that is $\parallel u \parallel_{W^{1,p}(V \times \D)} \le C$
 where a constant $C > 0$ is independent of $d$. Indeed, in the above estimates $\tau > 0$ and $\delta(\tau) > 0$ may be choosen ondependent of $d$.
 Since $u \in W^{1,p}(V \times \D)$ for any $d$ and $h(\bullet,d)$ is of class $C^1(\R^n)$ for every $d$, the composition $h(u,d)$ is of class $W^{1,p}(V \times \D$ as well and the chain rule can be applied for generalized derivatives: $h(u)_{x_j} = Dh(u) u_{x_j}$.

  Therefore we obtain 
 
 \begin{eqnarray*}
 \label{BishopEq9}
  \parallel  u - u^0  \parallel_{W^{1,p}(V \times \D)}  =  \parallel Th(u,d) \parallel_{W^{1,p}(V \times \D)}  \le    C   \parallel h(\bullet,d)  \parallel_{C^1_x}
  \end{eqnarray*}   
 where   $\parallel h(\bullet,d) \parallel_{C^1_x}$ denotes the $C^1$ norm of $h(x,d)$  with respect to the variable $x$, and $C > 0$ is a constant.
 But $\parallel h(\bullet,d)  \parallel_{C^1_x}  \to 0$ as $d  \to 0$,  and Lemma follows.

 \section{Proof of Theorem  \ref{MainTheo1}}

The following proposition is the key technical result.

\begin{prop}
\label{MainProp}
Let $M$ be an $m$-dimensional totally real manifold of class $C^1$ in  $\C^m$. Assume that $\Omega \subset \C^n$ is a pseudoconvex domain defined by $\Omega = \{ \phi < 0 \}$ where $\phi$ is a plurisubharmonic function of class $C^2$. and $d \phi \neq 0$ near $\Gamma = b\Omega$. Let also $W \subset \Omega$ 
 be a wedge (\ref{wedge})  with the edge $E \subset \Gamma$ of type (\ref{edge}).  Consider a holomorphic map   $f:  \Omega \to \C^m$ such that $f$ is continuos on $\Omega \cup E$ and $f(E) \subset M$.  Then,  for every $\delta > 0$,  and every $\alpha < 1$, the map $f$ extends to a H\"older 
$\alpha$-continuous mapping on $W_{\delta} \cup E$. 
\end{prop}
\proof Consider the family of discs constructed in the former section and attached to $E$ along $b\D^+$. It follows by the maximum principle that all discs belong to $\Omega$. 
Indeed, the flat discs fill a prescribed wedge of type (\ref{wedge}) with the edge $E_0 =  \R^n$.  More precisely, we can fix an open convex  cone $K$ in 
$W^0 = \{ (x,y)  \in  \R^{2n}: y_j < 0, j = 1,...,n  \}$ with the vertex at the origin and such that $ \overline K \cap r\B^n$ is contained in $W^0   \cup  \{ 0  \}$, for some $r > 0$ small 
enough. Clearly, the flat discs fill a neighborhood of $ \overline K  \cap r  \B^n$.  The same remains true for the cone $K_z$ obtained by the parallel translation of $K$ to the vertex at $z  \in  \R^n$. 
  Since the family $H(c,t,d)(\zeta)$ is a small perturbation of the flat discs in  $C^s(V \times \overline \D)$ (with any $0 < s < 1$), by continuity for $d $ small enough the family $H(c,t,d)$ also fills a presribed edge of type (\ref{shrinked}) with the edge $E_d$. By the holomorphic equivalence, the same is true for the initial edge $E_\delta$ with any $\delta > 0$.


 Applying the strong version of the Hopf lemma (see \cite{PiTs}) to the subharmonic function $\phi \circ H$ on $\D$, we 
obtain that $$\vert \phi \circ H(c,t)(\zeta) \vert \ge C (1 - \vert \zeta \vert)$$ with $C > 0$ independent of discs (i.e. on $(c,t)$) . Recall that 
$$C^{-1} dist (z, \Gamma) \le \vert \phi(z) \vert \le C dist(z,\Gamma)$$
Since $E \subset \Gamma$, we have 
$$dist(z, \Gamma) \le dist(z,E)$$

This implies the estimate

\begin{eqnarray}
\label{dist}
 1- \vert \zeta \vert \le C dist(H(c,t)(\zeta),E)
\end{eqnarray}

Recall that $M$ can be defined by  $M = \rho^{-1}(0)$, where $\rho$ is a nonegative strictly psh function of class $C^2$ (see \cite{Ch2,HW}). 
Now we apply the previous argument to each disc $f(z)$ with $z = H(c,t)(\zeta)$ and obtain

$$\rho(f(z))  \le C(1- \vert \zeta \vert) \le C dist(z,E)$$
Note that here the left estimate  is obtained in  \cite{Su1}.

Hence we obtain the key estimate

\begin{eqnarray}
\rho(f(z)) \le C dist(z,E)
\end{eqnarray}
for all $z \in W_\delta$.

With this estimate the argument  from \cite{Su1} (based on estimates of the Kobayashi metric in a tube neighborhood of $M$) literally goes through and gives that the mapping $f$ is
 $\alpha$-H\"older on $W_\delta$ for each $\alpha < 1$. This proves Proposition.

\bigskip 
\bigskip

Now Theorem \ref{MainTheo1} follows exactly as in \cite{Su1}.

\end{document}